\documentclass[11pt]{article}
\usepackage{latexsym,amssymb}
\usepackage[centertags]{amsmath}
\begin{document}

%%% Abbreviation %%%
\newcommand{\qed}{\hfill\Box}\newcommand{\bi}{\bigskip}
\newcommand{\sm}{\smallskip}
\newcommand{\w}{\sqrt}
\newcommand{\wh}{\widehat}
\newcommand{\wt}{\widetilde}
\newcommand{\ee}{\end{equation}}
\newcommand{\eea}{\end{eqnarray}}
\newcommand{\bean}{\begin{eqnarray*}}
\newcommand{\eean}{\end{eqnarray*}}
\newcommand{\suml}{\sum\limits}
\newcommand{\intl}{\int\limits}
\newif\ifpctex
\newcommand{\noi}{\noindent}
\newcommand{\vt}{\vartheta}
\newcommand{\vp}{\varphi}
\newcommand{\vr}{\varrho}
\newcommand{\ve}{\varepsilon}

%%%%%%%  Theorem environment  %%%%%%%
  \newtheorem{definition}{Definition}
  \newtheorem{prob}{Problem}
  \newtheorem{conjecture}{Conjecture}
  \newtheorem{quest}{Question}
  \newtheorem{theorem}{Theorem}
  \newtheorem{proposition}{Proposition}
  \newtheorem{lemma}{Lemma}
  \newtheorem{corollary}{Corollary}
  \newtheorem{assumption}{Assumption}
  \newtheorem{example}{Example}
  \newtheorem{remark}{Remark}
  \newtheorem{result}{Result}
%%%%%%%  Domains of numbers   %%%%%%%
  \newcommand{\R}{\mathbb{R}}
  \newcommand{\Z}{\mathbb{Z}}
  \newcommand{\N}{\mathbb{N}}
  \newcommand{\C}{\mathbb{C}}
  \newcommand{\CL}{\mathcal{L}}
  \newcommand{\CH}{\mathcal{H}}
  \newcommand{\CM}{\mathcal{M}}
  \newcommand{\CU}{\mathcal{U}}

\newcommand{\Rand}[1]{\marginpar{#1}}
     \renewcommand{\Rand}[1]{}
\marginparwidth2.5cm
\newcommand{\be}[1]{\Rand{\vspace{0,6cm}\tt #1}\begin{equation}\label{#1}}
\newcommand{\bea}[1]{\Rand{\vspace{0,6cm}\tt #1}\begin{eqnarray}\label{#1}}
\newcommand{\beL}[2]{\Rand{\vspace{0,6cm}\tt #1}\begin{lemma}[#2]\label{#1}}
\newcommand{\beA}[2]{\Rand{\vspace{0,6cm}\tt #1}\begin{assumption}[#2]\label{#1}}
\newcommand{\beD}[2]{\Rand{\vspace{0,6cm}\tt #1}\begin{definition}[#2]\label{#1}}
\newcommand{\beT}[2]{\Rand{\vspace{0,6cm}\tt #1}\begin{theorem}[#2]\label{#1}}
\newcommand{\beP}[2]{\Rand{\vspace{0,6cm}\tt #1}\begin{proposition}[#2]\label{#1}}
\newcommand{\beC}[2]{\Rand{\vspace{0,6cm}\tt #1}\begin{conjecture}[#2]\label{#1}}
\newcommand{\beCor}[2]{\Rand{\vspace{0,6cm}\tt #1}\begin{corollary}[#2]\label{#1}}

%%% Convergences %%%
  \newcommand{\atp}[2]{\genfrac{}{}{0pt}{}{#1}{#2}}
  \newcommand{\llto}{{_{\displaystyle{\atp{\ll}{t\to\infty}}}}}
  \newcommand{\Tto}{{_{\displaystyle{\atp{\Longrightarrow}{t\to\infty}}}}}
  \newcommand{\tto}{{_{\displaystyle{\atp{\longrightarrow}{t\to\infty}}}}}
  \newcommand{\tTo}{{_{\displaystyle{\atp{\longrightarrow}{T\to\infty}}}}}
  \newcommand{\Tno}{{_{\displaystyle{\atp{\Longrightarrow}{n\to\infty}}}}}
  \newcommand{\TNo}{{_{\displaystyle{\atp{\Longrightarrow}{N\to\infty}}}}}
  \newcommand{\tno}{{_{\displaystyle{\atp{\longrightarrow}{n\to\infty}}}}}
  \newcommand{\tMo}{{_{\displaystyle{\atp{\longrightarrow}{M\to\infty}}}}}
  \newcommand{\Mto}{{_{\displaystyle{\atp{\longrightarrow}{M\to\infty}}}}}
  \newcommand{\Nto}{{_{\displaystyle{\atp{\longrightarrow}{N\to\infty}}}}}
  \newcommand{\tNo}{{_{\displaystyle{\atp{\longrightarrow}{N\to\infty}}}}}
  \newcommand{\trhoo}{{_{\displaystyle{\atp{\longrightarrow}{\rho\to\infty}}}}}
  \newcommand{\Ttto}{{_{\displaystyle{\atp{\Longrightarrow}{T,t\to\infty}}}}}
  \newcommand{\Trhoo}{{_{\displaystyle{\atp{\Longrightarrow}{\rho\to\infty}}}}}

\title{\Large Law of Large Numbers for a Class of
Superdiffusions } \thispagestyle{empty}

\author{J{\'a}nos Engl\"ander$^1$, Anita Winter$^2$}
\date{}
\maketitle \hfill\today
\begin{abstract}
Pinsky (1996) \cite{Pin96} proved that  the finite mass
super\-diffusion $X$ corresponding to the semilinear operator
$Lu+\beta u-\alpha u^2$ exhibits local extinction  if and only if
$\lambda_c \le 0$, where $\lambda_c:=\lambda_c(L+\beta)$ is the
generalized principal eigenvalue of $L+\beta$ on $\R^d$. For the
case when $\lambda_c > 0$, it has been shown in Engl\"ander and
Turaev (2000) \cite{EngTur00} that in law the  super\-diffusion
locally behaves like $\exp[t\lambda_c]$ times a non-negative
non-degenerate random variable, provided that
 the operator $L+\beta-\lambda_c$ satisfies a certain spectral condition
(`product-criticality'), and that $\alpha$ and $\mu=X_0$ are `not
too large'.

In this article we will  prove that the  convergence in law used in
the formulation in \cite{EngTur00}  can actually be replaced by
convergence in probability. Furthermore, instead of  $\mathbb {R}^d$
we will consider a general Euclidean domain $D\subseteq \mathbb
{R}^d$.

%Our result generalizes the main theorem in \cite{FleSwa 03} and
%also proves a conjecture proposed in the same paper.

As far as the proof of our main theorem is concerned, the heavy
analytic method of \cite{EngTur00} is replaced by a different,
simpler and more probabilistic one. We introduce a space-time
weighted superprocess ($H$-transformed superprocess) and use it in
the proof  along with some elementary probabilistic arguments.

\vspace{1cm} \noindent
\end{abstract}
%\tableofcontents

\vspace*{1cm}

\footnoterule \noi \hspace*{0.3cm}{\footnotesize $^1$ { Department
of Statistics and Applied Probability, University of California,
Santa Barbara, CA 93106-3110, USA. \hspace*{0.3cm} Email:{\texttt{ englander@pstat.ucsb.edu}%
}\vspace{2pt}}\\
\hspace*{0.3cm}{\footnotesize $^2$ Mathematisches Institut,
Universit\"at Erlangen--N\"urnberg, Bismarckstra\ss{}e $1\frac12$,
91054 Erlangen, Germany. \hspace*{0.3cm} Email: {\tt
winter@mi.uni-erlangen.de}}}
\newpage
\hspace{4.3cm}{\bf R\'esum\'e:}

Pinsky (1996) \cite{Pin96} a prouv\'e que le processus de
superdiffusion de masse finie $X$ correspondant {\'a} l'operateur
semilin\'{e}aire $Lu+\beta u-\alpha u^2$ poss\`{e}de la
propri\'{e}t\'{e} d'extinction locale si, et seulement si,
$\lambda_c \le 0$, o\`{u} $\lambda_c:=\lambda_c(L+\beta)$ est la
valeure propre generalis\'{e}e de $L+\beta$ dans $\R^d$. Dans le
cas o\`{u} $\lambda_c > 0$, et pour un operatour
$L+\beta-\lambda_c$ possedant une condition spectrale (de
`criticalit\'{e}-produit'), et pourvu que  $\alpha$ et $\mu=X_0$
ne  soient pas trop grands, Engl\"ander and Turaev (2000)
\cite{EngTur00} ont montr\'e le processus se comporte localement
et en loi comme $\exp[t\lambda_c]$ avec une constante
multiplicative al\'eatoire non d\'eg\'en\'er\'ee.

Dans l'article pr\'esent, nous montrons que la convergence en loi
de \cite{EngTur00} peut \^{e}tre renforc\'ee par la convergence en
probabilit\'e. De plus, l'espace $\mathbb {R}^d$ est generalis\'e
a un domaine Euclidien quelconque.

S'agissant de la preuve du theor\`{e}me principale la lourde
methode analytique de \cite{EngTur00} est remplac\'ee par une
approche probabiliste plus simple. Nous introduisons une
renormalisation spatio-temporelle du superprocessus\linebreak
(`{\it $H$-transformed superprocess}') que nous utilisons dans la
preuve combin\'ee a des arguments probabilistes \'el\'ementaires.

\section{Introduction}
\subsection{Preparation}
We consider a superprocess which arises as the short life time and
high density diffusion limit of a branching particle system, which
can be described as follows: in the $n^{\mathrm{th}}$ approximation
step each particle has mass $1/n$ and lives a random time which is
exponential with mean $1/n$. While a particle is alive, its motion
is described by a diffusion process  corresponding to the operator
$L$. At the end of its life, the particle dies and is replaced by a
random number of particles situated at the parent particle's final
position. The distribution law of the number of descendants is
spatially varying such that the mean number of descendants is
$1+\frac{\beta(x)}{n}$, while the variance is assumed to be
$2\alpha(x)$. All these mechanisms are independent of each other.
The process is determined by the quadruple $(L,\beta,\alpha;D)$,
where $L$ is a second order elliptic operator corresponding to the
underlying diffusion process on $D$. See Appendix A in Engl\"ander
and Pinsky (1999) \cite{EngPin99} for a precise statement on the
particle approximation.

We start by presenting a formal description of the model
considered in this article. For convenience we first recall the
basic notation: let $D\subseteq \mathbb{R}^d$ be a domain and let
${\cal B}(D )$ denote the  Borel sets of $D$. We write ${\cal
M}_f(D)$ and ${\cal M}_c(D)$ for the class of finite measures
resp.\ the class of finite measures with compact support on ${\cal
B}(D )$. For $\mu\in{\cal M}_f(D )$, denote $ \|\mu\|:=\mu (D )$
and let $C^+_b(D )$ and
 $C^+_c(D )$ be the class of non-negative bounded continuous
resp.\ non-negative continuous  functions $D \rightarrow\R$ having
compact support. Write $C^{k,\eta}(D )$ for the usual H\"older
spaces of index $\eta\in (0,1]$ including derivatives of order $k$,
and set $C^{\eta}(D):=C^{0,\eta}(D)$.

We continue with the definition of the
$(L,\beta,\alpha;D)$-superdiffusion, $X$. Let $L$ be an elliptic
operator on the domain $D\subseteq\R^d$ of the form \be{A3}
   L:=\frac{1}{2}\nabla \cdot a \nabla +b\cdot\nabla,
\end{equation}
where $a_{i,j},b_i\in C^{1,\eta}(D)$, $i,j=1,...,d$, for some
$\eta\in (0,1]$, and the matrix $a(x):=(a_{i,j}(x))$ is symmetric,
and positive definite for all $x\in D$. In addition, let
$\alpha,\beta\in C^{\eta}(D)$, and assume that $\alpha$ is
positive, and $\beta$ is bounded from above.

\smallskip
We now present our model.
\begin{definition}[Time-homogeneous  superdiffusion]\qquad\qquad

Let $\left( X,\mathbf{P}^{\mu\,},\,\mu\in\mathcal{M}_f(D)\right) $
denote the $(L,\beta,\alpha;D )$-superdiffusion. That is, $X$ is
the unique $\mathcal{M}_f(D)$-valued continuous (time-homogeneous)
Markov process which satisfies, for any $g\in C^+_b(D)$ ,
\begin{equation}
   {\bf E}^{\mu}\exp\left\langle  X_{t\,},-g\right\rangle
 =\exp\, \langle\mu,-u(\cdot ,t)\rangle,
\label{Laplace.functional}%
\end{equation}
where $u$ is the minimal nonnegative solution to
\begin{equation}
\left.
\begin{array}
[c]{c}%
u_{t}=Lu+\beta u-\alpha u^{2}\quad\text{\emph{on}}\;D%
\times(0,\infty),\vspace{8pt}\\
\lim\limits_{t\downarrow 0}u(\cdot,t)=g(\cdot).
\end{array}
\;\right\}  \label{cum.equ}%
\end{equation}
As usual, $\left\langle \nu,f\right\rangle $ denotes the integral
$\int_{D}\nu(\mathrm{d} x)\,f(x).$
\end{definition}

(See Dynkin (1991, 2002) \cite{Dyn91}, \cite{Dy02} or Dawson
(1993) \cite{Daw93} for the definition of superprocesses in
general; see Engl\"ander and Pinsky (1999) \cite{EngPin99} for
more on the definition in the particular setting above.)

\begin{remark}[Time-inhomogeneous
superdiffusion] \rm
%Even though
The model under consideration is a time-homogeneous process.
However, it  is important to point out that for the formulation of
the main theorem and the proof, the introduction of certain {\it
time-inhomogeneous} superdiffusions is required. The previous
definition will therefore be generalized for time-inhomogeneous
superdiffusions in  Appendix B. (See Definition \ref{D2}.)
\end{remark}

Let \be{A5q}
\begin{aligned}
  \lambda _{c}
 &:=
  \lambda _{c}(L+\beta )
  \\
 &:= \inf \{\lambda \in \mathbb{R}\ : \
  \exists u>0\ \text{satisfying}\ (L+\beta -\lambda )u=0\ \text{in}\
  D \}
\end{aligned}
\end{equation}
denote the {\it generalized principal eigenvalue }for $L+\beta $
on $D$. Let $\xi^{L}$ be the diffusion process on $D$
corresponding to $L$, and denote by $\mathbb P^x$ the law
of $\xi^L$ starting at $x\in D$.
Then from a probabilistic point of view, the
generalized principal eigenvalue can be equivalently expressed as
\be{A5qq}
  \lambda_{c}
 =
  \sup_{\{A:\ A\subset \subset D,\ \partial A\ \mathrm{is}\ C^{2,\eta}\}}
  \lim_{t\to \infty }\frac{1}{t}\log{\mathbb{E}^{x}
  \Big[
  \exp\big[\int_{0}^{t}\beta (\xi^{L}_s)\,\mathrm{d}s\big];\,\tau^{A}>t
  \Big]},
\end{equation}
for any $x\in D,$ where $\tau ^{A}=\inf \{ t\geq 0:\xi^{L}(t)\not\in
A\}$, and the $C^{2,\eta}$-boundary is defined with the help of
$C^{2,\eta}$-maps in the usual way. (See Section 4.4 in Pinsky
(1995) \cite{Pin95} on the subject). Hence, since $\beta$ is bounded
from above, $\lambda_{c}<\infty$; and it is known from standard
theory that for any $\lambda \geq \lambda _{c}$, there exists a
function $0<f \in C^{2,\eta }(D)$ such that $(L+\beta )f =\lambda f
$ on $D$. (See Section 4.3 in Pinsky (1995) \cite{Pin95})

Pinsky proved that $X$ {\em exhibits local extinction}  (i.e., the
support of $X$ leaves any given bounded set, ${\bf P}^{\mu}$-a.s.
for each $\mu\in{\cal M}_c$) if and only if $\lambda_c \le 0$.
(See Theorem 6 and Remark 1 in Pinsky (1996) \cite{Pin96}.)
\smallskip

From now on we are interested in the situation where $X$ does not
exhibit local extinction. We therefore assume that $\lambda_c>0$.
We get a first rough impression about the local growth rate  by
the following statement  taken from Theorem 7(b) in Pinsky (1996)
\cite{Pin96}: \addtocounter{lemma}{-1} \begin{lemma}[Local
behavior in expectation]
\begin{itemize}
\item[{}] \hspace{2cm}

For $\mu\in{\cal M}_c(D)$, and $g\in C_c^+(D)$, satisfying
$||\mu||\neq 0$ and $g\not\equiv 0$, \be{A6}
   \limsup_{t\to\infty}\,\exp[-\rho t]{\bf E}^{\mu}[\langle X_t,g\rangle]
 =
   \left\{\begin{array}{cc}0&\mbox{ if }\rho >\lambda_c\\\infty&
   \mbox{ if }\rho <\lambda_c\end{array}\right..
\end{equation}
(See also Appendix \ref{Pinproblem}.)
\end{itemize}

\end{lemma}

\smallskip

We are therefore going to concentrate on scaling with the exponent
$\rho=\lambda_c$. In addition to  the concept of the generalized
principal eigenvalue we will then need some further ones  from the
so-called {\it criticality theory} of second order elliptic
operators. In particular, we will use the concepts of {\it critical}
and {\it product-critical} (or  product-$L^1$ critical) operators.
Recall that  the operator $L+\beta-\lambda_c$ is called critical if
there exists a positive function $f$ satisfying that
$(L+\beta-\lambda_c)f=0$ but there is no (minimal positive) Green's
function for the operator $L+\beta-\lambda_c$. In this case $f$ is
unique up to constant multiples and is called the {\em ground
state}. The operator $L+\beta-\lambda_c$ is called product-critical
if it is critical with  ground state $0<\phi_c$, and $\phi_c$ and
$\tilde{\phi_c}$ (i.e.\ the ground state for the formal adjoint of
$L+\beta-\lambda_c$) satisfy $\langle
\mathrm{d}x,\phi_c\tilde{\phi_c}\rangle<\infty$. In this case we
normalize them by $\langle \mathrm{d}x,\phi_c\tilde{\phi_c}\rangle
=1$.

If $L+\beta-\lambda_c$ possesses a Green's function, then it is
called {\it subcritical}.

For the reader, it will be handy to have Appendix 2 of
\cite{EngTur00} at hand, where a review on criticality theory is
given. For a complete presentation of the theory, the reader is
referred to Chapter 4 in \cite{Pin95}.\smallskip

\subsection{Motivation}

When $D=\mathbb{R}^d$ and  $L+\beta-\lambda_c$ is
product-critical, it is known from Theorem 1 in \cite{EngTur00}
that if $\parallel\alpha\phi_c\parallel_{\infty}<\infty$ and the
initial state $\mu$ is such that $\langle \mu,\phi_c\rangle
<\infty$, then  the following  holds in the vague topology:

\be{A7} \lim_{t\to\infty}  \exp[-\lambda_c t] X_t\, (\mathrm{d}x)
 =
    N^{\mu}\,\tilde{\phi_c}\, \mathrm{d}x\ \ \mathrm{in\ law},
\end{equation}
where the limiting non-negative non-degenerate random variable
$N^{\mu}$ was identified with the help of a certain {\em invariant
curve}.

It is important to point out that even though product-criticality
is equivalent to the ergodicity of an auxiliary diffusion process
(see next section), {\it the original motion process
corresponding to $L$
does not have to be ergodic}. In fact it can even be transient --
see Example 23 in \cite{EngTur00}.

(On the other hand, it follows from the discussion in Appendix A
that when $L+\beta-\lambda_c$ on $D\subseteq \mathbb R^d$ is not
product-critical, then for $\mu\in{\cal M}_c$, and $g\in
C_c^+(D)$, \be{A6o} \lim_{t\to\infty}\,\exp[-\lambda_c t]\langle
X_t,g\rangle
 = 0 \ \ \mathrm{in\ }L^1.)
\end{equation}

There are two disadvantages of the method used in \cite{EngTur00}.
First, the assumption that $D=\mathbb{R}^d$ is essential. Second,
the proof does not yield any probabilistic insight.

In this paper our goal is to improve the statement in (\ref{A7})
and to provide a proof that is probabilistic in nature. We will
show that the `{\it Law of Large Numbers}' holds, that is, that
one can replace the convergence in law by convergence in
probability\footnote{Since the limit is not constant, therefore,
unlike in classical probability theory, one has to distinguish
between convergence in law (WLLN) and convergence in probability
(LLN).}. Furthermore we will drop the assumption that
$D=\mathbb{R}^d$. In the proof we will replace the analytic
reasoning given in \cite{EngTur00} (which relies on dynamical
systems) by a more probabilistic one using space-time weighted
superprocesses ($H$-transforms).

We suspect that in fact the Strong Law of Large Numbers holds,
that is that the convergence in probability can  be replaced by
almost sure convergence. However we could not upgrade the proof of
this paper to give the Strong Law.

In the recent paper \cite{FleSwa 03} the authors study  a
supercritical superprocess taking values in the space of finite
measures on $[0,1]$, whose underlying motion is the Wright--Fisher
diffusion corresponding to the operator $$L=\frac{1}{2}
x(1-x)\frac{\mathrm{d}^2}{\mathrm{d}x^2},$$ and whose branching
mechanism is $\gamma u(1-u)$ with $\gamma
>0$ (that is, $\alpha=\beta=\gamma$). They establish a dichotomy in the
long-time behavior of this superprocess. For $\gamma\le 1$, the mass
in the interior $(0,1)$ dies out after a finite random time, while
for $\gamma >1$, the mass in $(0,1)$ grows exponentially with rate
$\gamma - 1$ (as $t\to\infty$ and with positive probability) and is
approximately uniformly distributed over $(0,1)$.

This result is in line with that of \cite{EngTur00} if one considers
the restriction of the superprocess on the (open) domain $(0,1)$. In
fact it is easy to show that $\lambda_c:=\gamma-1$ is the principal
eigenvalue of the linearized elliptic operator $L+\gamma$.
%$A+\gamma$ on $(0,1)$. I.e.
%$$\lambda_c (A+\gamma )=\lambda,$$ or equivalently,
%$$\lambda_c (A)= -1.$$
Here is a possible argument:  the operator $L+\gamma-\lambda_c=L+1$
can be $h$-transformed ($h=v$, where $v$ is an explicitly given
function in the paper) into a diffusion operator, which -- according
to their Lemma 20 -- corresponds to a (positive) recurrent
diffusion. Consequently, this $h$-transformed operator  is {\it
critical}, and thus  its principal eigenvalue is zero. By
$h$-transform invariance, the same is then true for the original
operator $L+\gamma-\lambda_c.$ (See again Chapter 4 in
\cite{Pin95}). Furthermore, the product-criticality and boundedness
assumptions are automatically satisfied by the boundedness of
$D=(0,1)$.

Finally the fact that the limiting measure is the Lebesgue measure,
is also in line with \cite{EngTur00}. Indeed, according to
\cite{EngTur00}, the limiting density is a harmonic function with
respect to the adjoint of $L+\gamma-\lambda_c=L+1$, that is with
respect to $\tilde L+1$, where $\tilde L$ is the adjoint of $L$. An
easy computation reveals that $$\tilde L+1=\frac{1}{2}
x(1-x)\frac{\mathrm{d}^2}{\mathrm{d}x^2}+(1-2x)\frac{\mathrm{d}}{\mathrm{d}x}.$$
Since the adjoint of a critical operator is also critical, and since
 positive harmonic functions for a critical
operator are unique up to constant multiples,  the limiting density
must be a properly normalized constant on the unit interval, that
is, the limiting density is one.

However, as the authors point out referring to \cite{EngTur00},
`\it their methods use in an essential way the fact that their
underlying space is $\mathbb{R}^d$ (and not an open subset of
$\mathbb{R}^d$, like $(0, 1)$), and therefore their results are
not applicable to our situation.'\rm

In the present article, as already mentioned,  we manage to
overcome this difficulty, so the result of \cite{FleSwa 03} will
fit our main result (in \cite{FleSwa 03} the result is somewhat
stronger as they prove convergence in $L^2$).

%In \cite{FleSwa 03} the authors propose as a conjecture (see their
%formula (23)) that the scaling limit they obtained in $L^2$ is in
%fact an almost sure limit. This will too follow from our main
%result.

\section{Main Result}
Recall from (\ref{A5q}) the definition of the principal eigenvalue
$\lambda_c$  of $L+\beta$ on $D$ and the corresponding ground state
$\phi_c$, and that throughout the paper we assume that
$\lambda_c>0$. Also, $\{{\cal S}_s\}_{s\ge 0}$ will denote the
semigroup (`expectation semigroup') corresponding to the operator
$L+\beta$ on $D$. So far we have recalled (\ref{A7}). In order to
replace in (\ref{A7}) the convergence in law by convergence in
probability, we will assume the same conditions as in
\cite{EngTur00}, except that we work with a generic domain $D$.

\begin{assumption}\label{assum}
\begin{itemize}
\item[{}] \hspace{2cm}

In addition to the assumption that $\lambda_c>0$, also assume that
$L+\beta-\lambda_c$ is product-critical, that $\alpha\phi_c$ is
bounded and that $X$ starts in a state $\mu$ with
$\langle\mu,\phi_c\rangle<\infty$.
\end{itemize}
\end{assumption}
Before   reading the remainder of this section, it is recommended
that the reader consults Appendix B regarding the definition of
time-inhomogeneous superdiffusions as well as the space-time
$H$-transform (weighted superdiffusion) introduced there.

Let $X$ be a $(L,\beta,\alpha;D)$-superdiffusion with $X_0=\mu$.
 Let $H(x,t):=\exp (-\lambda_c t )\phi_c(x),\ x\in D,\
t\ge 0$. It turns out (see Lemma~\ref{LH1} in Appendix B) that the
(time-inhomogeneous) process $X^H$ defined by \be{H1aa}
  X_t^H:=H(\cdot,t)X_t\quad \left(\mbox{that is,}\
  \frac{\mathrm{d}X_t^H}{\mathrm{d}X_t}=H(\cdot,t)\right), \quad
  t\geq 0
\end{equation}
is an $(L+a\frac{\nabla \phi_c}{\phi_c}\cdot \nabla,0,\alpha \phi_c
e^{-\lambda_c t};D)$-superdiffusion. In the sequel $\widetilde{\bf
E}$ and $\widetilde{\mathrm{Var}}$ will denote expectation and
variance with respect to the law of $X^H$.

\beL{JA}{Bounded variance}
\begin{itemize}

 \item[{}] \hspace{.2cm} \be{JL(A2)i}
\lim_{t\to\infty}\widetilde{\mathrm{Var}}^{\phi_c\mu}
(\|X_t^H\|)=\int_0^{\infty}\mathrm{d}s\,e^{-2\lambda_c s}\,
 \langle \mu,{\cal S}_s[\alpha \phi_c^2]\rangle<\infty,\end{equation}
and  $\| X^H\|$ is a uniformly integrable $\widetilde{\bf
P}^{\phi_c\mu}$-martingale.
\end{itemize}
\end{lemma}
\bf{Proof:}\rm\ Let $\overline{X}^H$ denote the {\em total mass
process}, i.e., \be{Y11}
   \overline{X}^H:=\| X^H\|.
\end{equation}   Abbreviate \be{LH}
   L^{\phi_c}_0:=L+a\frac{\nabla \phi_c}{\phi_c}\cdot \nabla
\end{equation}
and note that in fact
$$L^{\phi_c}_0(u)=\phi_c^{-1}(L+\beta-\lambda_c)(\phi_c u)=
   H^{-1}(L+\beta+\partial_t)(Hu).$$
(Here $\partial_t$ denotes differentiation with respect to time.)
Let ${\cal S}^{\phi_c}$ denote the $h$-transformed semigroup with
$h=\phi_c$, that is ${\cal S}^{\phi_c}_s(\cdot)=(\phi_c)^{-1}{\cal
S}_s(\phi_c \cdot).$

Define ${\cal S}^H_s:=
   e^{-\lambda_c s}{\cal S}^{\phi_c}_s$;
then the semigroup $\{{\cal S}^H_s\}_{s\ge 0}$ corresponds to the
operator $L_0^{\phi_c}$ that has no zeroth order part. In
particular then
\begin{equation}\label{contraction} {\cal S}^H_s
1\le 1.
\end{equation}
Finally, the product-criticality assumption on $L+\beta-\lambda_c$
guarantees that the diffusion process corresponding to
$L^{\phi_c}_0$ on $D$ is positive recurrent (ergodic) (see Section
4.4. in Pinsky \cite{Pin95}). (Since ergodicity implies
conservativeness, thus in fact ${\cal S}^H_s 1= 1;$ nonetheless, for
us it will be enough to know (\ref{contraction}).)

By Lemma~\ref{LH1}(b) of Appendix B along with  Theorem A2 in
\cite{EngPin99}, we have that for all $f\in
C^2_{\mathrm{const}}(D):=\{f\in C^2(D):\ \exists
\Omega\subset\subset D\mathrm{\ such\ that}\ f=\mathrm{const}\
\mathrm{on}\ D\setminus \Omega \}$, \be{H3a} \mathrm{d}\langle
X^H_t,f\rangle =\langle
X^H_t,L^{\phi_c}_0f\rangle\,\mathrm{d}t+\mathrm{d}M_t(f),
\end{equation}
where $\{M_t(f)\}_{t\ge 0}$ is a square-integrable $\widetilde{{\bf
P}}^{\phi_c\mu}$-martingale,
%com%{\bf APPLICABLE TO TIME DEPENDENT SETUP?}
and its quadratic variation (i.e. the increasing process in the
Doob-Meyer decomposition) $\langle M(f)\rangle$ is given by \be{H4}
   \langle M(f)\rangle_t
 =
   \int^t_0\mathrm{d}s\,e^{-\lambda_c s}\langle X^H_s,\alpha
   {\phi_c}f^2\rangle,\ \ t\ge 0.
\end{equation}
(One can take the function class $C^2_{\mathrm{const}}(D)$ instead
of just $C^2_c(D)$, because the diffusion process corresponding to
$L^{\phi_c}_0$ on $D$ is conservative, that is, it never leaves
the domain $D$ with probability one.)

Applying (\ref{H3a}) to the function $f\equiv 1$, it follows that
$\overline{X}^H$ is a $\widetilde{{\bf P}}^{\phi_c\mu}$-martingale.
Furthermore,  by (\ref{H4}), \be{H4f<}
%A%
  \widetilde{{\bf E}}^{\phi_c\mu}\left[\langle X^H_t,1\rangle^2\right]
 =
  \langle\mu, \phi_c\rangle^2+\int_0^{t}\mathrm{d}s\,e^{-\lambda_c s}\,
   \langle \phi_c\mu,{\cal S}^H_s[\alpha \phi_c]\rangle.
\end{equation}
That is \be{H4f<2} \widetilde{\mathrm{Var}}^{\phi_c\mu} (\|X_t^H\|)
 =\int_0^{t}\mathrm{d}s\,e^{-\lambda_c s}\,
   \langle \phi_c\mu,{\cal S}^H_s[\alpha \phi_c]\rangle=
   \int_0^{t}\mathrm{d}s\,e^{-2\lambda_c s}\,
 \langle \mu,{\cal S}_s[\alpha \phi_c^2]\rangle.
\end{equation}
Letting $t\to\infty$ we obtain the first statement of the lemma.

Replacing $t$ by $\infty$ in the first of the integrals in
(\ref{H4f<2}), we have from (\ref{contraction}) and from our
 assumptions that
$$\widetilde{\mathrm{Var}}^{\phi_c\mu} (\|X_t^H\|)\le\int_0^{\infty}\mathrm{d}s\,e^{-\lambda_c s}\,
   \langle \phi_c\mu,{\cal S}^H_s[\alpha \phi_c]\rangle
\le{\lambda_c }^{-1}\parallel\alpha \phi_c\parallel_{\infty}\,
   \langle \mu, \phi_c\rangle<\infty.$$
Hence, by (\ref{H4f<}),
%A%
$$\sup\nolimits_{t\ge 0}\widetilde{{\bf E}}^{\mu\phi_c}\left[\langle
X^H_t,1\rangle^2\right]<\infty,$$ and consequently $\overline{X}^H$
is uniformly integrable. This completes the proof of the second
statement of the lemma. $\hfill\square$

\begin{remark}\rm\
 Our proof of LLN will indeed use the condition that
 $\alpha \phi_c$ is  bounded, however it is quite possible
 that this condition is not necessary and that assuming the finiteness of the integral in
 (\ref{JL(A2)i}) (along with $\langle\mu,
\phi_c\rangle<\infty$) would suffice.
\end{remark}

An immediate consequence of uniform integrability is that
%A%
$\widetilde{{\bf
E}}^{\phi_c\mu}[\overline{X}^H_\infty]=\langle\mu,\phi_c\rangle$,
which is finite by assumption, and positive for $\mu\not =0$. This
yields that $\widetilde{{\bf
P}}^{\phi_c\mu}[\overline{X}^H_\infty=0]<1$ for $\mu\not =0$. We
record this in a lemma.
 \beL{LH1q}{Limit of the total mass}
\begin{itemize}
\item[{}]\hspace{2cm}

The martingale $\overline{X}^H$ has a $\widetilde{{\bf
P}}^{\phi_c\mu}$-a.s. limit
$\overline{X}^H_{\infty}:=\lim_{t\to\infty}\overline{X}^H_t$ which
is positive with positive probability.
\end{itemize}
\end{lemma}
\subsection{Heuristics for the Law of Large Numbers}
Before stating the {\it Law of Large Numbers} for the class of
superdiffusions under consideration,  in this subsection we give
some heuristic computations. These will justify why we call our main
result `the Law of Large Numbers'.

Working with the $H$-transformed superprocess and at the same time,
having the particle approximation in mind, consider particles with
underlying motion $Y$ corresponding to the elliptic operator
$L_0^{\phi_c}$ (the  probabilities for
$Y$ will be denoted by
$\{\mathbb P^{x}, \ x\in D\}$) and with critical binary branching at
rate $\exp [-\lambda_c t]\alpha(x)$ at position $x\in D$ and time
$t\ge 0$. Furthermore let the system be started with initial
discrete measure being ``close'' to $\nu\neq 0$.

Let ${\cal I}_t^n$ denote the collection of particles alive at
time $t$ in the $n^{\mathrm{th}}$ approximation step.
%, and for each $1\le k\le |{\cal
%I}_t^n|$, let $x_t^k$ denote the position of the
%``$k$th'' among them.
Finally, the event of `survival' is $\{|{\cal I}^n_t|>0\ \forall
t>0\}$.

The Law of Large Numbers would mean that if $0\not\equiv f\in
C_c^+(D)$, then as $t\to\infty$, (and without further specifying
what ``$\approx$'' means), \be{AA1}
 \|\nu\| \frac{\frac{1}{|{\cal I}^n_t|}\sum_{x\in{\cal I}_t^n}f(x)}
  {\langle \nu,{\mathbb E}^{x}[f(Y_t)]\rangle}\approx 1\hspace{0.5cm}
 \hspace{0.5cm}  \mathrm{on}\ \ \{|{\cal I}^n_t|>0\ \forall t>0\}.
\end{equation}
 Now recalling that in the $n^{th}$
approximating step the individual particle mass is scaled down by
$n$ and recalling also Lemma~\ref{LH1q}, one has that (for $n$
large), $|{\cal I}^n_t|\approx n\overline{X}^H_{\infty}$ as
$t\to\infty$. Putting this together with (\ref{AA1}), one gets
formally, that for large $n$, \be{AA1b}
 \|\nu\| \frac{\frac{1}{n}\sum_{x\in{\cal I}_t^n}f(x)}
  {\langle \nu,{\mathbb E}^{x}[f(Y_t)]\rangle}\approx \overline{X}^H_{\infty},\ \
 \ \mathrm{as}\ t\to\infty.
\end{equation}
Note that in fact
$$\langle \nu,{\mathbb
E}^{x}[f(Y_t)]\rangle=\widetilde{\mathbf E}^{\nu}\langle
X^H_t,f\rangle=e^{-\lambda_c t}\,\mathbf E^{\mu}\langle
X_t,f\phi_c\rangle$$ ($\nu=\phi_c\mu$). (The first equality can be
shown for instance by taking first the particular case
$\nu=\delta_x,\ x\in D$, and using that the two expectations satisfy
the same parabolic problem; then integrating with respect to
$\nu(\mathrm{d}x)$.) Furthermore, passing to the limit (as
$n\to\infty$) formally, the numerator of the fraction on the left
hand side of (\ref{AA1b}) becomes
$$\langle X^H_t,f\rangle=e^{-\lambda_c t}\langle X_t,f\phi_c\rangle.$$
Hence, for the new test function $0\not\equiv\hat f:=f\phi_c\in
C_c^+(D)$,
$$\frac{ \langle X_t,\hat f\rangle}{\mathbf{E}^{\mu}\langle X_t,\hat
f\rangle}\approx\frac{\overline{X}^H_{\infty}}{\|\nu\|}.$$

\subsection{Main theorem}
Making  the  intuition of the previous subsection precise, we now
state our main result:
\beT{JC1}{Law of Large Numbers}
\label{LLNproduct-critical}
\begin{itemize}
\item[{}] \hspace{2cm}

Let $f\in C_c^+(D)$. If $f\not\equiv 0$ and $||\mu||\neq 0$, then
\be{A71} \lim_{t\to\infty}\, \frac{\langle X_t,f\rangle}{{\bf
E}^{\mu}\langle X_t,f\rangle}
=\frac{\overline{X}^H_{\infty}}{\langle\mu,\phi_c\rangle}, \qquad
\mathrm{in }\;\,{\bf P^{\mu}}\mathrm{-probability}.
\end{equation}
\end{itemize}
\end{theorem}\smallskip

Comparing our theorem with (\ref{A7}), we can now identify the
limiting distribution: $N^{^\mu}=\overline{X}^H_{\infty}$ in law.
\begin{remark}\rm
\begin{itemize}
\item[{}] \hspace{2cm}

One has to be a bit careful though when making heuristic
inferences  using the particle picture as in the previous
subsection.

Obviously, the {\it discrete} system in the $n^{th}$ approximation
step is so that $\liminf_{t\to\infty}|{\cal I}^n_t|\ge 1$ under
survival. That is, $\liminf_{t\to\infty}\frac{1}{n}|{\cal I}^n_t|>0$
under survival. Recall that, heuristically,  (for $n$ large),
$\frac{1}{n}|{\cal I}^n_t|\approx \overline{X}^H_{\infty}$ as
$t\to\infty$.

On the other hand, in the recent paper \cite{Eng04} an example of a
{\it superprocess} is given that satisfies the conditions of our
previous theorem and for which
$${\bf P}^{\mu}(\overline{X}^H_{\infty}=0\mid S)>0,$$
where $S$ is the event of survival, $S:=\{\|X_t\|>0,\ \forall t\ge
0 \}.\hfill\Diamond$
\end{itemize}
\end{remark}
\begin{conjecture}\rm
\begin{itemize}
\item[{}] \hspace{2cm}

We conjecture that convergence in probability in
 (\ref{A71}) can be replaced by almost sure
convergence.
\end{itemize}
\end{conjecture}
We close this section with a remark concerning an old result.
\begin{remark}\rm

\begin{itemize}
\item[{}] \hspace{2cm}

A simple case of a superdiffusion is when $D=\mathbb R^d,\ d\ge
1,\ L=\frac{1}{2}\Delta$, with $\alpha,\beta$ positive constants
(supercritical super-Brownian motion). Here $\lambda_c=\beta$ and
$$\frac{1}{2}\Delta+\beta-\lambda_c= \frac{1}{2}\Delta.$$ Since
$\phi_c=\tilde \phi_c\equiv 1,\,d\ge 1,$ the operator
$\frac{1}{2}\Delta$ is either critical but not product-critical
($d\le 2$), or subcritical ($d\ge 3$). Therefore this case is not
included in our setup. On the other hand, the corresponding
(Strong) Law of large Numbers is well known -- at least for the
discrete particle systems. Using techniques from Fourier transform
theory, Watanabe (\cite{Wat67}) proved SLLN for branching-Brownian
motion in $\mathbb R^d$ and in certain subdomains of it. It is not
clear however if his method can be generalized for more general
branching diffusions. $\hfill\diamond$
\end{itemize}
\end{remark}
\section{Proof of the result}
The  proof is  based on two observations. The first one is that
the problem can be formulated in terms of $X^H$, that is, one can
reduce the problem to the investigation of a {\it critical}
superdiffusion with {\it ergodic motion component} and {\it
exponentially decaying branching rate} (again, recall that $X^H$
is an $(L+a\frac{\nabla \phi_c}{\phi_c}\cdot \nabla,0,\alpha
\phi_c e^{-\lambda_c t};D)$-superdiffusion).

The second one is that by considering some large time $t+T$ (where
both $t$ and $T$ are large), the changes of $\overline{X}^H$ are
negligible after time $t$, while  the remaining time  $T$ is still
long enough to distribute the produced mass according to the ergodic
flow given by the $H$-transformed migration.

To simplify  notation, we will write $W:=X^H$ (and, accordingly,
$\overline {W}_{\infty}:=\overline {X}^H_{\infty}$). Denote
$\nu:=W_0$. By assumption, $\|\nu\|=\langle \mu
,\phi_c\rangle<\infty.$ We need to show that for all $\epsilon>0$,
\begin{equation}\label{W}
\lim_{t\to\infty}{\bf P}^{\mu}\Big(\Big|\frac{\langle
X_t,f\rangle} {{\bf E}^{\mu}\langle
X_t,f\rangle}-\frac{\overline{W}_{\infty}}{\|\nu\|}\Big|>\epsilon\Big)=0.
\end{equation}
%With a slight abuse of notation, from now, until the end of the
%proof, $\bf P$ will denote the probability corresponding to the
%$H$-transformed process $W$ rather than the original process $X$.
%(Recall that $\|\nu\|<\infty$.)
Recall that $\widetilde{{\bf P}}$ denotes the probabilities with respect
to the law of $W$.
Denoting $\phi_c^{-1}f=:g\in
C_c^+(D)$, we rewrite (\ref{W}) as
%A%
$$\lim_{t\to\infty}\widetilde{{\bf P}}^{\nu}\Big(\Big|\frac{\langle W_t,g\rangle}
{\widetilde{{\bf E}}^{\nu}\langle
W_t,g\rangle}-\frac{\overline{W}_{\infty}}{\|\nu\|}\Big|>\epsilon\Big)=0.$$
It is easy to check that the limiting invariant density for the
diffusion corresponding to $L_0^{\phi_c}$ is $\phi_c\tilde\phi_c$
(recall the normalization
$\int_D\mathrm{d}x\,\phi_c(x)\tilde\phi_c(x)=1$). Since $Z:=\widetilde{{\bf E}}^{\nu}W$ is just the deterministic
$L_0^{\phi_c}$-flow starting from $\nu$, therefore
$$\lim_{t\to\infty}\widetilde{{\bf E}}^{\nu}\langle
W_t,g\rangle=||\nu||\cdot\langle \phi_c\tilde\phi_c ,g \rangle,$$
and consequently our statement is tantamount to saying that for
all $\epsilon>0$,
$$\lim_{t\to\infty}\widetilde{{\bf P}}^{\nu}\left(\left|\langle W_t,g\rangle-
\langle \phi_c\tilde\phi_c ,g
\rangle\overline{W}_{\infty}\right|>\epsilon\right)=0.$$  Let
$T>0$ and let $Z_{W_t}$ denote the deterministic
$L_0^{\phi_c}$-flow starting from the (random) measure $W_t$. Then
\begin{equation}\label{three}
\widetilde{{\bf P}}^{\nu}\big(\big|\langle W_{t+T},g\rangle-\langle
\phi_c\tilde\phi_c ,g
\rangle\overline{W}_{\infty}\big|>\epsilon\big)\le
S_1(t)+S_2(t,T)+S_3(t,T),
\end{equation} where
\begin{align}
S_1(t)&:=\widetilde{{\bf P}}^{\nu}\Big(\big|\langle \phi_c\tilde\phi_c ,g
\rangle\overline{W}_{\infty}-\langle \phi_c\tilde\phi_c ,g
\rangle\|W_{t}\|\big|>\epsilon/3\Big)\nonumber\\
S_2(t,T)&:=\widetilde{{\bf P}}^{\nu}\Big(\big|\langle \phi_c\tilde\phi_c ,g
\rangle\|W_{t}\|-\langle
Z_{W_t}(T),g\rangle\big|>\epsilon/3\Big)\nonumber\\
S_3(t,T)&:=\widetilde{{\bf P}}^{\nu}\big(\big|\langle
Z_{W_t}(T),g\rangle-\langle
W_{t+T},g\rangle\big|>\epsilon/3\big).\nonumber
\end{align}
Take $\limsup_{t\to\infty}\limsup_{T\to\infty}$ on both sides of
(\ref{three}).  We have
$$\limsup_{t\to\infty}\widetilde{{\bf
P}}^{\nu}\left(\left|\langle W_{t},g\rangle -\langle
\phi_c\tilde\phi_c ,g
\rangle\overline{W}_{\infty}\right|>\epsilon\right)\le I+II+III,$$
where
\begin{align}
I&:= \limsup_{t\to\infty}S_1(t)\\
II&:=\limsup_{t\to\infty}\limsup_{T\to\infty}S_2(t,T)\\
III&:=\limsup_{t\to\infty}\limsup_{T\to\infty}S_3(t,T).
\end{align}
Now, since $\|W_t\|\to\overline{W}_{\infty}$ as $t\to\infty$
a.s.,
$$I=\lim_{t\to\infty}S_1(t)=0.$$
Also $ II=0,$ because for all fixed $t\ge 0,\
\lim_{T\to\infty}S_2(t,T)=0.$ (Indeed, for {\it all}
$\omega\in\Omega$, $\lim_{T\to\infty} \langle
Z_{W_t(\omega)}(T),g\rangle=\langle \phi_c\tilde\phi_c ,g
\rangle\|W_{t}(\omega)\|$.) Therefore, if we show that
\begin{equation}\label{S3}
III=0,
\end{equation} then we are done.

In order to do this, use  at time $t$ that $W$ is a
time-inhomogeneous Markov-process, and then apply Chebysev's
inequality:
\begin{equation}
\label{MarkovChebysev}
\begin{aligned}
   S_3(t,T)
 &=
   \widetilde{{\bf E}}^{\nu}{\bf \hat P}^{W_t}\left(\left|\langle Z_{W_t}(T),
   g\rangle-\langle W_{t+T},g\rangle\right|>\epsilon/3\right)
  \\
 &\le
   9 \epsilon^{-2}\,\widetilde{{\bf E}}^{\nu} \mathrm{\hat{\sigma}^2}_{W_t}
   \langle W_{T},g\rangle,
\end{aligned}
\end{equation}
where ${\bf \hat P}$ is the law of the $(L+a\frac{\nabla
\phi_c}{\phi_c}\cdot \nabla,0,\alpha \phi_c e^{-\lambda_c
(t+s)};D)$-superdiffusion (here $t$ is fixed and $s>0$ is time)
and $\hat{\sigma}^2$ denotes variance.

Let us now recall  how the formulae for the first two moments of
$\langle W_{T},g\rangle$ are obtained: by writing
$u_{\theta}(t,x)$ for the solution of the semilinear parabolic
evolution equation (corresponding to the superprocess) with
initial value $\theta g$, one differentiates (repeatedly) with
respect to $\theta$ and sets $\theta=0$.

For  time-homogeneous processes with constant branching rate  this
is written down in detail in \cite{Etheridge 02}, p.37. Since the
derivation of these `moment formulae' only uses differentiation with
respect to $\theta$ (and not $t$ or $x$), therefore the proof goes
through for the more general setting where coefficients are
space-time-dependent.

In our case, from these moment formulae and from
(\ref{MarkovChebysev}), one obtains (recall also
(\ref{contraction})) that for all $T>0$,
\begin{equation}
\label{moments}
\begin{aligned}
   S_3(t,T)
 &\le
   \frac{9}{\epsilon^{2}}\cdot\widetilde{\mathbf{E}}^{\nu}\int_0^T\mathrm{d}s\,
   2 e^{-\lambda_c(t+s)}
   \big\langle W_t,\mathcal{S}^H_s\big[\alpha\phi_c
   \big(\mathcal{S}^H_{T-s}g\big)^2\big]\big\rangle
  \\
 &\le
   \frac{C\cdot\widetilde{{\bf E}}^{\nu}\|W_t\|}
   {\epsilon^{2}\cdot \lambda_c e^{\lambda_c t}},
\end{aligned}
\end{equation}
 with
$$C=C(\|g\|,\|\alpha\phi_c\|):=18\,\|
\alpha\phi_c\|\cdot\|g\|^2.$$ (Note that we have an extra factor $2$
relative to \cite{Etheridge 02} --- indeed in \cite{Etheridge 02}
the nonlinear term in the semilinear parabolic evolution equation is
$\frac{1}{2}\gamma u^2$.) Recall that $\|W\|$ is a $\widetilde{{\bf
P}}^{\nu}$-martingale with mean $\|\nu\|$ and continue
(\ref{moments}) with
$$=\frac{C\cdot\|\nu\|}{\epsilon^{2}\cdot \lambda_c e^{\lambda_c t}}.$$
 Since this holds for
all $T> 0$, thus $$\limsup_{T\to\infty}S_3(t,T)\le
\frac{C\cdot\|\nu\|}{\epsilon^{2}\cdot \lambda_c e^{\lambda_c
t}}.$$ Letting $t\to\infty$, one obtains (\ref{S3}), completing
our proof. $\qquad\qquad\square$

\begin{appendix}
\section{The behavior of the process in expectation}
\label{Pinproblem}\rm For the cases when  $L+\beta-\lambda_c$ is
subcritical, or critical but $\langle
\mathrm{d}x,\phi_c\tilde{\phi_c}\rangle=\infty$, Theorem 7(b)(ii)
in \cite{Pin96} states that for $g\in C_c^{+}$, \be{A8}
 \lim_{t\to\infty}   e^{-\lambda_c t}{\bf E}^{\mu}[\langle X_t,g\rangle]=0.
\end{equation}
Note, however, that in the proof there is a glitch: the proof
simply refers to Theorem 4.9.9 in \cite{Pin95} which deals with
the product-critical case only, and is therefore not applicable
for the cases mentioned.

Nevertheless,  for the subcritical case, and for
$\mu\in\mathcal{M}_c(D)$, the statement can be verified by a very
simple argument as follows (cf. Theorem 4.9.1. in \cite{Pin95}).
First, note that by the first moment formula,
$$e^{-\lambda_c t}{\bf E}^{\mu}[\langle
X_t,g\rangle]=\langle\mu ,[\mathcal{T}_t g](x)\rangle,$$  where
$\{\mathcal{ T}_t\}_{t\ge 0}=\{e^{-\lambda_c t}\mathcal{
S}_t\}_{t\ge 0}$ denotes the semigroup corresponding to the
operator $L+\beta-\lambda_c$ on $D$.

Make an $h$-transform now: $$[\mathcal{T}_t
g](x)=h(x)[\mathcal{T}^h_t (\hat g)](x),$$ where $\hat {g}:=g
h^{-1}$, and pick an $h>0$ satisfying $(L+\beta-\lambda_c)h=0$, to
reduce the problem to the proof of \be{A8a} \lim_{t\to\infty}
\langle  h\mu, \mathbb E^{x} [\hat g(\xi^h_t)]\rangle= 0,
\end{equation}
 where $\xi^h$ denotes the diffusion corresponding to
$\mathcal{T}^h$, that is, to the operator $L_{0}^{h}$ (defined
analogously to $L_{0}^{\phi_c}$ in (\ref{LH}) with $\phi_c$
replaced by $h$) and $\mathbb  E$ denotes the corresponding
expectation. (Of course, $h\mu\in\mathcal{M}_c(D)$).

Furthermore, it is enough to show the statement with
$\mu:=\delta_x,\ x\in D$, because once we know that,  the general
statement follows by bounded convergence: $\mathbb E^{x} [\hat
g(\xi^h_t)]\le \|\hat g\|$ for all $x\in D$ and $t\ge 0$.

In \cite{Pin95}, Chapter 4 it is shown that subcriticality is
invariant under $h$-transforms, and that the transience of a
diffusion is equivalent to the subcriticality of the corresponding
elliptic operator. Therefore, in our case, it follows that $\xi^h$
is transient. Since $g$ is compactly supported, by an obvious
comparison argument, it is enough to show (\ref{A8a}) with $g$
replaced by the indicator $1_{B}$, where $B$ is a ball. By
transience $1_{B}(\xi^h_t)\to 0\ \mathrm{as}\ t\to\infty$ a.s.,
and the statement follows from this and bounded convergence.

Similarly, the critical but non-product critical case can be reduced
to the analogous (but much subtler) problem of showing (\ref{A8a})
for a null-recurrent $\xi^h$. (Cf. the well known analogous limit
theorem for countable state space Markov chains -- see e.g.
Proposition 5.3 and Corollary 6.39 in \cite{KemSneKna76}.) In fact,
as mentioned in the notes at the end of Chapter 4 in \cite{Pin95},
this result is known in the case when $L$ is  {\it symmetric}   with
respect to some reference measure $\rho\mathrm{d}x$ (see
\cite{ChK91} or \cite{Sim93}). (Recall that $L$ is symmetric if and
only if $b=a\nabla Q$ for some $Q\in C^{2,\eta}(D)$, $\eta\in
(0,1]$, and in this case $L$ possesses a self-adjoint extension due
to the Friedrichs extension theorem -- see Chapter 4 in
\cite{Pin95}.) Recently Pinchover completed the result by  proving
it  for the general (non-selfadjoint) setting (see \cite{Pinch04}).

Consequently,   the $\rho >\lambda_c$ (over-scaling) part of
(\ref{A6}) is  immediate. One does not need however the above deep
result for the  $\rho >\lambda_c$ part. Here is a simple
alternative proof:  using an h-transform with an $h>0$ satisfying
$(L+\beta-\lambda_c)h=0$, the statement is equivalent to \be{A8b}
\lim_{t\to\infty} e^{(\lambda_c -\rho)t}\mathbb
E^{x}[g(\xi^h_t)]=0,
\end{equation}
for each $x\in D$, which is true in virtue of the boundedness of
$g$.

The $\rho <\lambda_c$ (under-scaling) part of (\ref{A6}) is
harder, and we are only able to provide the rigorous proof of the
somewhat weaker assertion:
\be{A6b}
   \limsup_{t\to\infty}\,\exp[-\rho t]{\bf E}^{\mu}[\langle X_t,g\rangle]
 = \infty.
\end{equation}

To this end, denote by $p^{L+\beta}(t;x,\cdot),\ x\in D,$ the
kernel corresponding to the operator $L+\beta$ and note that since
$g$ is compactly supported, by an obvious comparison argument, it
is enough to prove that \be{A8d}
   \limsup_{t\to\infty} \frac{1}{t}\log \int_B \mu(\mathrm{d}x)\, p^{L+\beta}(t;x,B)
 =
   \lambda_{c},
\end{equation}
for each $x\in D$, and Borel set $B\subset\subset D$. Clearly, we
may assume that $\|\mu\|=1$. To verify (\ref{A8d}),  make again an
$h$-transform with an $h>0$ satisfying $(L+\beta-\lambda_c) h=0$ on
$D$. Then $L+\beta$ transforms into $L_0^h+\lambda_c$. Moreover,
since the generalized principal eigenvalue is invariant under
$h$-transforms and since
$\lambda_c(L_0^h+\lambda_c)=\lambda_c(L_0^h)+\lambda_c$, it follows
that $\lambda_c(L_0^h)=0$. Let  $p_0^h(t,x,\cdot)$ denote the
transition measures corresponding to $L_0^h$. Fix $x\in D$ and
$B\subset\subset D$. Since $h( \cdot
)h^{-1}(x)p^{L+\beta}(t,x,\cdot)=e^{\lambda_c t}p_0^h(t,x,\cdot)$
(see Theorem 4.1.1. in \cite{Pin95}), and since $h$ is bounded
between two positive constants on $B$, we have \be{A8e}
   \limsup_{t\to\infty} \frac{1}{t}\log\int_B \mu(\mathrm{d}x)\, p_0^h(t;x,B)
 =
   \limsup_{t\to\infty} \frac{1}{t}\log\int_B \mu(\mathrm{d}x)\, p^{L+\beta}(t;x,B)-\lambda_c.
\end{equation}
Since $\|\mu\|=1$ and $p_0^h(t;x,B)\le 1$, the left hand side of
(\ref{A8e}) is non-positive, giving immediately \be{A8f}
   \lambda_{+}(x,B)
 :=
   \limsup_{t\to\infty} \frac{1}{t}\log\int_B \mu(\mathrm{d}x)\, p^{L+\beta}(t;x,B)\le \lambda_c.
\end{equation}
 Suppose now that $\lambda_{+}(x,B)< \lambda_c$ and pick
\be{A8g}
   c
 \in (-\lambda_c,-\lambda_{+}(x,B)).
\end{equation}
Then by (\ref{A8e}) and (\ref{A8g}), along with Fubini's theorem,
one obtains \be{A8h1} \int_B \mu(\mathrm{d}x)\,\int_0^{\infty} e^{c
t}p^{L+\beta}(t,x,B)\, \mathrm{d}t<\infty.
\end{equation}
Hence, for almost every $x\in B$, \be{A8h}  G^{(L+\beta+c)}(x,B)
 :=\int_0^{\infty} e^{c t}p^{L+\beta}(t,x,B)\,
\mathrm{d}t<\infty.
\end{equation}
It follows from general theory then  that
$e^{ct}p^{L+\beta}(t,x',B')$ is in fact integrable for {\it  all }
$x'\in D$ and $B'\subset\subset D$, that is, that  the operator
$L+\beta+c$ possesses a (minimal positive) {\it Green's function} on
$D$; however
 this contradicts to the  well known fact that an  operator  with positive generalized
principal eigenvalue does not  possess a Green's function. (In our
case $\lambda_c(L+\beta+c)=\lambda_c+c>0.$)

\section{The $H$-transform of superdiffusions}
\label{Htransformsection}\rm This section treats a generalization of
the $h$-transform for superdiffusions introduced in \cite{EngPin99}.
The $h$-transform was used in the proofs in \cite{EngTur00}. The
method used in the present paper however requires  the spatial
function $0<h=h(x)$ to  be replaced by a space-time function
$0<H=H(x,t)$. (The  reader should not confuse with the space-time
harmonic transformation yielding a Girsanov-type change of measure
-- see  e.g. \cite{Ove94}.)

We start with the more general definition of a  time-inhomogeneous
superdiffusion. Let $\tilde{L}$ be an elliptic operator on
$D\times\R^+$ of the form \be{A3n}
   \tilde{L}:=\frac{1}{2}\nabla \cdot \tilde{a} \nabla +\tilde{b}\cdot\nabla
\end{equation}
where the functions
$\tilde{a}_{i,j},\tilde{b}_i:D\times\R^+\to\R$, $i,j=1,...,d$ are
$C^{1,\eta}(D)$  (for some $\eta\in (0,1]$) in the space, and
continuously differentiable in the time coordinate. Moreover
assume that the symmetric matrix $\tilde{a}(x,t):=(a_{i,j}(x,t))$
is  positive definite for all $x\in D$ and $t\in\R^+$.

In addition, let $\tilde{\alpha},\tilde{\beta}:D\times\R^+\to\R$,
be $C^{\eta}(D)$ in the space, and continuously differentiable in
the time coordinate. Finally assume that $\tilde{\alpha}$ is
positive, and $\tilde{\beta}$ is bounded from above.

\beD{D2}{Time-inhomogeneous
$(\tilde{L},\tilde{\beta},\tilde{\alpha};D)$-superdiffusion} \rm\
\begin{itemize}

\item[{}] {\bf (i)} The
{\rm$(\tilde{L},\tilde{\beta},\tilde{\alpha};D)$-superdiffusion}
is a  measure-valued (inhomogeneous) Markov process, $(X,{\bf
P}^{\mu,r};\,\mu\in{\cal M}_f(D),r\ge 0)$,  that is, a family
$\{{\bf P}^{\mu,r}\}$ of probability measures where ${\bf
P}^{\mu,r}$ is a probability on $C([r,\infty) )$ and the family is
indexed by ${\cal M}_f(D))\times [0,\infty)$, such that the
following holds: for each $g\in C_b^+(D)$ and $\mu\in{\cal
M}_f(D)$, \be{A4c}
   {\bf E}^{\mu,r}[\exp-\langle X_t,g\rangle]=\exp-\langle\mu,u(\cdot,r;t,g)\rangle,
\end{equation}
where $u=u(\cdot,\cdot;t,g)$ is a particular  non-negative solution
to the backward equation \be{A5s}
\begin{aligned}
   -\partial_r u&=\tilde{L}u+\tilde{\beta} u-\tilde{\alpha} u^2\hspace{1cm}
  \mbox{in }D\times (0,t),
  \\
   \lim_{r\uparrow t}u(\cdot,r;t,g)&=g(\cdot).
\end{aligned}
\end{equation}
{\bf (ii)} To determine the solution $u$ uniquely, use the
equivalent {\it forward} equation along with the minimality of the
solution: fix $t>0$ and introduce the `time-reversed' operator
$\hat{L}$ on $D\times (0,t)$ by \be{A3n2}
   \hat{L}:=\frac{1}{2}\nabla \cdot \hat{a} \nabla
   +\hat{b}\cdot\nabla,
\end{equation}
where, for $r\in[0,t]$,
$$\hat a(\cdot,r):=\tilde a(\cdot,t-r)\ \mathrm{and}\ \hat b(\cdot,r):=\tilde b(\cdot,t-r);$$
furthermore let
$$\hat \beta(\cdot,r):=\tilde \beta(\cdot,t-r)\ \mathrm{and}\ \hat \alpha(\cdot,r):=
\tilde \alpha(\cdot,t-r).$$ Consider now $v$, the {\it minimal}
non-negative solution to the {\it forward} equation\be{A5s2}
\begin{aligned}
   \partial_r v&=\hat{L}v+\hat{\beta} v-\hat{\alpha} v^2\hspace{1cm}
  \mbox{in }D\times (0,t),
  \\
   \lim_{r\downarrow 0}v(\cdot,r;t,g)&=g(\cdot).
\end{aligned}
\end{equation}
Then $$u(\cdot,r;t,g)=v(\cdot,t-r;t,g).$$ (See also the remark
following this definition.)
\end{itemize}
\end{definition}
\begin{remark}[Minimal non-negative solutions for forward equations]\begin{itemize}

\item[{}]\rm
In the time-homogeneous case, minimal non-negative solutions for
forward equations have been constructed in Appendix
A in \cite{EngPin99}, and in Section 2 in \cite{EngPin03}
--- the construction uses the approximation of $D$ by compactly embedded subdomains
with Dirichlet condition on their boundaries. The construction goes
through for the time-inhomogeneous setting.

In \cite{EngPin99}
 and \cite{EngPin03} the time interval is $[0,\infty)$ rather than
 $[0,t]$. However that does not make any difference --- in fact the
 solution on the infinite time interval can be defined by first
 working on finite time intervals and then showing that they can be
 taken arbitrarily large without having the solution blown up.
 \end{itemize}
 \end{remark}

As we will see in Lemma \ref{LH1} (b), one way of defining a
time-inhomogeneous superdiffusion is to start with a
time-homogeneous one, and then to apply an `{\it $H$-transform}'. In
general, the $H$-transform of a time-inhomogeneous superdiffusion is
defined as follows. Let $0<H\in C^{2,\eta}(D)\times
C^{1,\eta}(\mathbb R^+)$ and let $X$ be a
$(\tilde{L},\tilde{\beta},\tilde{\alpha};D)$-superdiffusion. We
define a new process $X^H$ by
\begin{equation}\label{newprocess}
  X_t^H:=H(\cdot,t)X_t\quad \left(\mbox{that is,}\
  \frac{\mathrm{d}X_t^H}{\mathrm{d}X_t}=H(\cdot,t)\right), \quad
  t\geq 0.
\end{equation}
This way one obtains a new superdiffusion, which, in general, {\it
is not finite measure-valued} but only $\sigma$-finite
measure-valued. That is, if $\mathcal{M}(D)$ denotes the family of
all (finite or infinite) measures on $D$, then
$$X^H_t\in\mathcal{M}_H^{(t)}(D):=\{\nu\in \mathcal{M}(D)\mid H(\cdot,t)^{-1}\nu\in
\mathcal{M}_f(D)\}$$ (c.f. \cite{EngPin99}, p. 688.)

In \cite{EngPin99}, Section 2, it was shown, that, from an
analytical point of view, the (spatial) $h$-transform of the
superdiffusion is given by a certain transformation of the
corresponding semilinear operator. This remains the case for the
space-time $H$-transform.
\beL{LH1}{H-transform}
\begin{itemize}
\item[{}] \hspace{2cm}

Let $X^H$ be defined by (\ref{newprocess}). Then \item[(a)] $X^H$
is a $(\tilde{L}+\tilde{a}\frac{\nabla H}{H}\cdot\nabla,
  \tilde{\beta}+\frac{\tilde{L}H}{H}+\frac{\partial_rH}{H},
  \tilde{\alpha}H;D)$-superdiffusion.
\item[(b)] In particular, if $X$ is a time-homogeneous
  $(L,\beta,\alpha;D)$-superdiffusion, and $H$ is of the form
\be{H}
   H(x,t):=e^{-\lambda_c t}h(x),
\end{equation}
where $\lambda_c$ is the principal eigenvalue of $L+\beta$, and
$h$ is a positive solution of $(L+\beta)h=\lambda_c h$, then $X^H$
is a $(L+a\frac{\nabla h}{h}\cdot \nabla,0,\alpha h e^{-\lambda_c
t};D)$-superdiffusion.
\end{itemize}
\end{lemma}
\begin{remark}[Unbounded $\tilde \beta$'s]\begin{itemize}
\item[{}] \hspace{2cm}\rm

As it is already the case with the spatial $h$-transform for
superdiffusions, it is possible that the coefficient $\tilde
\beta$ transforms into a new coefficient that is no longer
bounded. In fact this can be the very definition of
superdiffusions for certain unbounded $\tilde \beta$'s (see
\cite{EngPin99}, Section 2 for explanation).

\end{itemize}\end{remark}

 {\bf Proof of
Lemma~\ref{LH1}. } In order to avoid minor technical
inconveniences, we implement the method in \cite{EngPin99} (see
the second paragraph on p. 689.). Namely, we use that the Laplace
transition functional restricted to the family of measures ${\cal
M}_c(D)$ and the family of functions $C^+_c(D)$ uniquely
determines a measure-valued Markov process, and we choose working
with these smaller spaces rather than replacing ${\cal M}_f(D)$
and $C^+_b(D)$ by $H$-dependent spaces.

Pick $\nu\in{\cal M}_c(D)$, and $f\in C_c(D)$. Define
$\mu^{(s)}:=\nu/H(\cdot,s)\in{\cal M}_c(D)$, and
$g^{(t)}(\cdot):=H(\cdot,t)f(\cdot)\in C^+_c(D)$. Obviously,
 \be{H3b}
  \widetilde{{\bf E}}^{\nu,s}\left[\exp-\langle X^H_t,f\rangle\right]
  ={\bf E}^{\mu^{(s)},s}\left[\exp-\langle
  X_t,g^{(t)}\rangle\right].
\end{equation}
By the log-Laplace equation (\ref {A4c}), we can continue with
$$=\exp-\left\langle\frac{\nu}{H(\cdot,s)},u(\cdot,s;t,g^{(t)})
\right\rangle.$$ Consider the operator  $${\cal A}:
C^{2,\eta}(D)\times C^{1,\eta}(\mathbb R^+)\mapsto C^{\eta}(D)\times
C^{\eta}(\mathbb R^+)$$ defined by \be{H4a}
\begin{aligned}
   {\cal A}(u)&:=\partial_s
   u+(\tilde{L}+\tilde{\beta})u-\tilde{\alpha}u^2.
\end{aligned}
\end{equation}
Define the $H$-transformed operator  ${\cal A}^H$ in the usual
 way:
\begin{equation}\label{conjugate}
 {\cal A}^H(u):=
   \frac{1}{H}{\cal A}(Hu).
\end{equation}
    Then a direct computation gives
\be{H5}
\begin{aligned}
   {\cal A}^H(u):=
   \frac{\partial_sH}{H}u+\partial_su+\tilde{L}u+\tilde{a}\frac{\nabla H}{H}\cdot
   \nabla u+\tilde{\beta}u+\frac{\tilde{L}H}{H}u-\tilde{\alpha}Hu^2.
\end{aligned}
\end{equation}
Another, trivial computation yields that if $$
v(\cdot,\cdot;t,f):=u(\cdot,\cdot;t,H(\cdot,t)f)/H(\cdot,t),$$ then
$v(\cdot,\cdot;t,f)$ is the solution in (\ref{A5s}) with ${\cal A}$
replaced by ${\cal A}^H$, and with the property in Definition
\ref{D2}(ii).
 Thus the quadruple
$(\tilde{L},\tilde{\beta},\tilde{\alpha};D)$ transforms into the
quadruple given in part (a).\sm

\medskip
\noi Part (b) is straightforward computation. $\qed$
\begin{remark}\rm
\begin{itemize}
\item[{}] \hspace{2cm}

It is precisely equation (\ref{conjugate}) that justifies the name
`$H$-transform'; the transformation on the semilinear operator works
the same way as Doob's $h$-transform would work on a linear
operator.
\end{itemize}
\end{remark}
\end{appendix}

\noi {\bf Acknowledgment. }\quad
The authors are grateful to the anonymous referee for
several comments and suggestions
that helped to improve the presentation and correctness
of the paper.


\begin{thebibliography}{40}
\bibitem{ChK91}{\bf Chavel, I., Karp, L. (1991)} Large time behavior
  of the heat kernel: the parabolic $\lambda$-potential alternative,
  Comment. Math. Helv., {\bf 66}, 541-556.
%\bibitem{Cro1957}{\bf Croft, H. T. (1957)} A question of limits, Eureka, {\bf 20}, 11-13.

\bibitem{Daw93}
{\bf  Dawson, D. A. (1993)}
\newblock Measure-valued {M}arkov processes.
\newblock In P.L. Hennequin, editor, {\em {\'E}cole d'\'et\'e de probabilit\'es
  de Saint Flour XXI--1991}, volume 1541 of {\em Lecture Notes in Mathematics},
  pages 1--260. Springer-Verlag, Berlin.

\bibitem{Dyn91}
{\bf  Dynkin, E. B.. (1991)}\newblock Branching particle systems and
superprocesses, \newblock {Ann. Probab.}, {\bf 19(3)}, 1157--1194.

\bibitem{Dy02}
{\bf Dynkin, E. B. (2002)}  \newblock{Diffusions, superdiffusions
and partial differential equations}, American Mathematical Society
Colloquium Publications, \bf 50\rm .
 American Mathematical Society, Providence, RI.

\bibitem{Eng04}{\bf Engl\"ander, J. (2004)}
An example and a conjecture concerning scaling limits of
superdiffusions,  Statist. Probab. Lett. {\bf 66(3)},  363-368.

\bibitem{EngPin99}{\bf Engl\"ander, J. and Pinsky, R.G. (1999)}
  On the construction and  support properties of measure-valued
  diffusions on $D\subset\mathbb R^d$ with spatially dependent branching,
  Ann. Probab., {\bf 27(2)}, 684-730.
\bibitem{EngPin03}{\bf Engl\"ander, J. and Pinsky, R.G. (2003)}
 Uniqueness/nonuniqueness for positive solutions to semilinear equations
 of the form $u_t=Lu+Vu-\gamma u^P$ in $R^n$,
 J. Diff. Equations, {\bf 192(2)}, 396-428.

\bibitem{EngTur00}{\bf Engl\"ander, J. and Turaev, D. (2002)} A scaling limit
  theorem for a class of superdiffusions,   Ann. Probab. {\bf 30(2)}, 683--722.
%\bibitem{Hur64}{\bf Hurwitz, A. (1964)} Vorlesungen \"uber allgemeine
%  Funktionentheorie. Springer-Verlag.




\bibitem{Etheridge 02}  {\bf Etheridge, A. (2000)}
\newblock{ An introduction to superprocesses.}
University Lecture Series, \bf 20.\rm\ AMS, Providence, RI.


\bibitem{FleSwa 03}{\bf  Fleischmann, K. and  Swart, J. (2003)}  Extinction
versus exponential growth in a supercritical super-Wright-Fisher
diffusion.  Stochastic Process. Appl. \bf 106(1)\rm, 141--165.

 \bibitem{KemSneKna76}{\bf Kemeny, J.G., Snell, J.L., and Knapp, A.W. (1976)} Denumerable Markov chains,
  Springer-Verlag.

%\bibitem{King63}{\bf Kingman, J. F. C. (1963)} Ergodic properties of
%continuous-time Markov processes and their discrete skeletons.
%Proc. London Math. Soc. {\bf 3(13)}, 593--604.


\bibitem{Ove94}{\bf Overbeck, L. (1994)} Some aspects of the Martin boundary
  of measure-valued diffusions. In Measure-valued processes, stochastic
  partial differential equations, and interacting systems, 179-186, AMS,
  Providence, RI.

\bibitem{Pinch04}{\bf Pinchover,\ Y. (2004)} Large time behavior of the heat
kernel, J. Functional Analysis {\bf 206(1)}, 191--209.

\bibitem{Pin95}{\bf Pinsky,\ R.G. (1995)} Positive harmonic functions and diffusion.
  Cambridge University Press.
\bibitem{Pin96}{\bf Pinsky,\ R.G. (1996)} Transience, recurrence and
  local extinction properties of the support for supercritical finite
  measure-valued diffusions, Ann. Probab., {\bf 24(1)},
  237-267.
\bibitem{Sim93}{\bf Simon, B. (1993)} Large time behavior of the heat
  kernel: On a theorem of Chavel and Karp, Proc. Amer. Math. Soc., {\bf
  118(2)}, 513-514.
\bibitem{Wat67}  {\bf Watanabe, S. (1967)} Limit Theorems for a Class of
Branching Processes. In Markov Processes and Potential Theory (ed.
J. Chover),  205--232. Wiley, New York.
\end{thebibliography}
\end{document}